\documentclass[10pt]{amsart}
\usepackage{amsmath, latexsym}
\usepackage{amssymb,amsthm}
\usepackage{amscd, 
}
 
\newtheorem{theorem}{Theorem}[section]
\newtheorem{proposition}[theorem]{Proposition}
\newtheorem{lemma}[theorem]{Lemma}

\theoremstyle{remark}
\newtheorem{remark}[theorem]{Remark}

\theoremstyle{definition}

\def\C{\mathbb{C}}
\def\real{\mathbb{R}}
\def\integer{\mathbb{Z}}
\def\supp{\mathrm{supp}}
\def\sphere{\mathbf{S}^{d-1}}
\def\cone{\mathbf{C}} 
\def\BB{\mathcal{B}}
\def\FF{\mathcal{F}}
\def\QQ{\mathcal{Q}}
\def\SS{\mathcal{S}}

\newcommand{\sob}[1]{_{W_*^{#1}}}%
\newcommand{\hol}[1]{_{C_*^{#1}}}%
\def\L{\mathcal{L}}
\def\P{\mathcal{P}}
\def\g{g}

\newcommand{\cout}[1]{} 

\begin{document}

\title[Anisotropic H\"older and Sobolev spaces] 
{Anisotropic H\"older and Sobolev spaces 
for hyperbolic diffeomorphisms}

\author{Viviane Baladi and Masato Tsujii}

\address{CNRS-UMR 7586, Institut de Math\'e\-ma\-ti\-ques  Jussieu, Paris, France}
\email{baladi@math.jussieu.fr}

\address{Department of Mathematics, Hokkaido University, 
Sapporo, Hokkaido, Japan}

\email{tsujii@math.sci.hokudai.ac.jp}

\date{January 2006, revised version}

\begin{abstract}
We study spectral properties of transfer operators for diffeomorphisms 
$T:X\to X$ on a Riemannian manifold $X$: 
Suppose that $\Omega$ is an isolated hyperbolic 
subset for $T$, with a compact isolating neighborhood  $V\subset X$. We first introduce 
Banach spaces of distributions supported on $V$, which are anisotropic versions of 
the usual space of $C^p$ functions $C^p(V)$ and of the generalized Sobolev spaces $W^{p,t}(V)$, 
respectively.  Then we show that the transfer operators associated
to $T$ and a smooth weight $g$ extend boundedly to these spaces, and 
we give bounds on the essential spectral radii of such extensions in terms of  hyperbolicity 
exponents. 
\end{abstract}

\thanks{VB thanks Scuola Normale Superiore Pisa
for hospitality, and Artur Avila  for useful
comments. Special thanks to Gerhard
Keller and S\'ebastien Gou\" ezel whose remarks helped
us to correct mistakes in a previous
version. MT thanks NCTS(Taiwan) for hospitality during his stay.}

\maketitle


\section{Introduction}
\label{intro}

Let $X$ be a $d$-dimensional $C^\infty$ Riemannian manifold and let $T:X \to X$ be a  
diffeomorphism which is of class $C^1$ at least. For a given complex-valued 
continuous function $g$ on 
$X$, we define the Ruelle transfer operator $\L_{T,g}$ by
\[
\L_{T,g}:C^0(X)\to C^0(X),\quad
\L_{T,g} u(x)=g(x) \cdot u\circ T(x).
\]
Such operators appear naturally in the study of fine statistical properties of 
dynamical systems and provide efficient methods, for instance, to estimate of decay 
of correlations. (We refer e.g. to \cite{Ba0}.)
Typical examples are the pull-back operator 
\begin{equation}\label{pullback}
T^* u:=\L_{T,1}u=u\circ T ,
\end{equation}
and the Perron-Frobenius operator
\begin{equation}\label{PFoperator}
\P u:=\L_{T^{-1},|\det DT^{-1}|}u=|\det DT^{-1}|\cdot u\circ T^{-1}.
\end{equation}
This paper is about  spectral properties of the operator $\L_{T,g}$. 

We shall require a hyperbolicity assumption on the mapping~$T$: 
Let $\Omega\subset X$ be a compact isolated invariant subset for $T$, with a compact isolating 
neighborhood $V$, that is,  $\Omega=\cap_{m\in \mathbb{Z}}T^{m}(V)$. We assume that $\Omega$ 
is a hyperbolic subset, that is, there exists an invariant decomposition 
$T_{\Omega}M=E^u\oplus E^s$ of the tangent bundle over $\Omega$,
satisfying 
$\|DT^m|_{E^s}\|\le C\lambda^{m}$ 
and 
$\|DT^{-m}|_{E^u}\|\le C\lambda^m$, for all $m\ge 0$ and $x\in \Omega$, with  constants 
$C>0$ and $0<\lambda<1$. Up to decomposing $\Omega$,
we may suppose that the dimensions of $E^u(x)$ and $E^s(x)$ are constant.

We define two local hyperbolicity exponents for each  $x\in \Omega$ and each $m\ge 1$ by 
\begin{equation}\label{hypexp}
\begin{aligned}
\lambda_{x}(T^{m})&=\sup_{v\in E^s(x)\setminus\{0\}} 
\frac{\|DT^{m}_x(v)\|}{\|v\|}\le C\lambda^m \quad \mbox{and}\\
\nu_x(T^{m})
&=\inf_{v\in E^u(x)\setminus\{0\}} 
\frac{\|DT^m_x (v)\|}{\|v\|}\ge C^{-1}\lambda^{-m}.
\end{aligned}
\end{equation}

Let $\omega$ be the Riemannian volume form on $X$, and let $|\det DT|$ be the Jacobian of 
$T$, that is, the function given by $T^*\omega=|\det DT|\cdot \omega$. Put, for $m\ge 1$, 
\[
\g^{(m)}(x)=
\prod_{k=0}^{m-1} g(T^k(x)).
\]
For real numbers $q\le 0\le p$, and for $1\le t\le \infty$, we set
\[
R^{p,q,t}(T, \g, \Omega, m)= 
 \sup_{\Omega}
 |\det DT^{m}_x|^{-1/t} |\g^{(m)}(x)|
\max\bigl \{ (\lambda_{x} (T^m))^{p}, (\nu_{x} (T^m))^{q}
\bigr \}, 
\]
where we read $(\cdot)^{1/\infty}=1$ for $t=\infty$. 
As $\log R^{p,q,t}(T, \g, \Omega, m)$ is sub-additive with respect to $m$, we have
\[
R^{p,q,t}(T,\g,\Omega):=\lim_{m\to \infty}\root{m}\of{R^{p,q,t}(T, \g,\Omega,m)}
=\inf_{m\ge 1}\root{m}\of{R^{p,q,t}(T, \g,\Omega,m)}.
\]

In this paper, we introduce  Banach spaces of distributions supported on~$V$,
show that the transfer operators $\L_{T,g}$ extend boundedly to 
those spaces and then give bounds for the essential spectral radii of these transfer
operators, using the quantities $R^{p,q,t}(T,\g,\Omega)$ introduced above. 
The main feature in our approach is that 
we work in Fourier coordinates. 
The definition and basic properties of the Banach spaces will be given in later sections. 
Here we state the main theorem as follows. 
Let $r_{ess} (L|_{\BB})$ be the essential spectral radius
of a bounded linear operator $L:\BB\to \BB$. 
For non-integer $s>0$, a mapping is of class $C^s$
if all its partial derivatives of order $[s]$ are 
$(s-[s])$-H\"older.

\begin{theorem}\label{main}
Suppose that $T$ is a $C^r$ diffeomorphism for a real number $r>1$, and let $\Omega$ be a
hyperbolic invariant set with compact isolating neighborhood $V$, as described above. 
Then, for any real numbers $q<0<p$ with $p-q<r-1$, there exist Banach spaces
$C^{p,q}_*(T,V)$ and $W^{p,q,t}_*(T,V)$, for $1<t<\infty$, of  distributions supported on $V$, 
such that, for any $C^{r-1}$ function $g:X\to \C$ supported on $V$, 
\begin{enumerate}
\item (H\"older spaces)
$\L_{T,g}$ extends boundedly  to $\L_{T,g}:C^{p,q}_*(T,V)\circlearrowleft$ and
\[
r_{ess}(\L_{T,g}|_{C^{p,q}_*(T,V)})\le 
R^{p,q,\infty}(T,\g,\Omega).
\]
\item (Sobolev spaces)
 $\L_{T,g}$ extends boundedly to 
$\L_{T,g}:W^{p,q,t}_*(T, V)\circlearrowleft$ and 
\[
r_{ess}(\L_{T,g}|_{W^{p,q,t}(T,V)}) 
\le  R^{p,q,t}(T,\g,\Omega) \, , \forall \, 1<t<\infty.
\]
\end{enumerate}
The Banach spaces $C^{p,q}_*(T,V)$ and $W^{p,q,t}_*(T,V)$ contain the set $C^s(V)$ of $C^s$ 
functions supported on $V$ for any $s>p$. 
\end{theorem}

\smallskip

{\bf Hyperbolic attractors and SRB measures.}
Let us see how to apply Theorem \ref{main} to a hyperbolic attractor:
Assume in addition to the above that $\Omega$ is an attracting hyperbolic
set and take the isolating neighborhood $V$ so that 
$T(V)\subset \mbox{interior}(V)$. Consider the pull-back operator $T^*$ defined by (\ref{pullback}) and the Perron-Frobenius operator $\P$ defined by (\ref{PFoperator}).
Note that these operators 
are adjoint to each other:
\begin{equation}\label{adjoint}
\int_{X} T^* u\cdot v \; d\omega=\int_{X} u \cdot \P v \;d\omega.
\end{equation}
Let $h:X\to [0,1]$ be a $C^\infty$ function supported on $V$ and satisfying $h\equiv 1$ 
on $T(V)$. 
Then the action of the 
operator $\mathcal{L}_{T^{-1},g}$ with $g(x)=|\det DT^{-1}(x)| \cdot h(x)$ coincides with that 
of the Perron-Frobenius operator $\P$ on $C^{p,q}_*(T^{-1},V)$ and 
$W^{p,q,t}_*(T^{-1},V)$, $1<t<\infty$. Therefore Theorem \ref{main} easily implies:

\begin{theorem}\label{PF} Let $\Omega$ be a hyperbolic attractor for a $C^r$ diffeomorphism 
$T:X\to X$ with $r>1$, and let $V$ be a compact neighborhood of $\Omega$ such that 
$T(V)\subset \mbox{{\rm interior}} (V)$ and  $\cap_{m\ge 0}T^m(V)=\Omega$. For real numbers $q<0<p$ with $p-q<r-1$, the Perron-Frobenius operator 
$\P$ extends boundedly to $\P:C^{p,q}_*(T^{-1},V)\to C^{p,q}_*(T^{-1},V)$ and also to 
$\P:W^{p,q,t}_*(T^{-1},V)\to W^{p,q,t}_*(T^{-1},V)$, and it holds
\[
r_{ess}(\P|_{C^{p,q}_*(T^{-1},V)})\le R^{-q,-p,1}(T,1,\Omega),
\]
and
\[
r_{ess}(\P|_{W^{p,q,t}_*(T^{-1},V)})\le R^{-q,-p,t/(t-1)}(T,1,\Omega)
\quad \mbox{for $1<t<\infty$.}
\]
(The above bound is strictly less than $1$ if $t>1$ is close enough to
$1$.)

For real numbers $q<0<p$ with $p-q<r-1$, the 
modified pull-back operator $T^*_h u:= h\cdot (u\circ T)$ 
extends boundedly to 
$T^*_h:C^{p,q}_*(T,V)\to C^{p,q}_*(T,V)$ and also to 
$T^*_h:W^{p,q,t}_*(T,V)\to W^{p,q,t}_*(T,V)$, 
and it holds
\[
r_{ess}(T_h^*|_{C^{p,q}_*(T,V)})\le R^{p,q,\infty}(T,1,V)<1 .
\]
Also,
\[
r_{ess}(T^*_h|_{W^{p,q,t}_*(T,V)})\le R^{p,q,t}(T,1,V), 
\quad \mbox{for $1<t<\infty$.}
\]
(The above bound is strictly less than $1$ if $t$ is large
enough.)
\end{theorem}

Once we have the estimates in Theorem \ref{PF}, 
it is not difficult to see that the spectral radius 
of the modified pull-back operator 
$T^*_h$ on  $C^{p,q}_*(T,V)$, and on $W^{p,q,t}_*(T,V)$ 
for large enough $t$, are equal to 
one ($h$ is a fixed point of $T^*_h$). 
If $(T,\Omega)$ is topologically mixing in addition, then $1$ is the 
unique eigenvalue on the unit circle, it is a simple eigenvalue, and
the fixed vector of the dual of
$T^*_h$ gives rise to  the SRB measure $\mu$ on 
$\Omega$: This 
corresponds to exponential decay of correlations for $C^p$ observables 
and $\mu$. (See Blank--Keller--Liverani \cite[\S3.2]{BKL} for example.)

\begin{remark}
{}From (\ref{adjoint}) and Theorem \ref{PF}, if $T$ is Anosov, the 
Perron-Frobenius operator 
$\P$ acts naturally on the (strong) dual spaces of $C^{p,q}_*(T,X)$ and $W^{p,q,t}_*(T,X)$, 
and we have for instance
\[
r_{ess}(\P|_{(C^{p,q}_*(T,X))^*})\le R^{p,q,\infty}(T,1,X)<1 ,
\]
for real numbers $q<0<p$ with $p-q<r-1$.
\end{remark}

\smallskip

{\bf Spectral stability.}
We  point out that, in the setting of Theorem \ref{main}, 
there is $\epsilon >0$  so that
if $\widetilde T$ and $\widetilde g$, respectively, are
$\epsilon$-close to $T$ and $g$, respectively, in the $C^r$, resp.
$C^{r-1}$, topology,  then the associated
operator $\L_{\widetilde T,\widetilde g}$
has same spectral properties than $\L_{T,g}$
on {\it the same Banach spaces.} Spectral
stability can then be proved, as it has been done \cite{BKL} or \cite{GL}
for the norms defined there (see also the historical comments below).

\smallskip 

{\bf Organization of the paper.}
After defining 
a version of our norms in $\real^d$
in Section \ref{S1}, we proceed in the usual way: prove compact
embeddings in Section \ref{S2} and a Lasota-Yorke
type estimate in Section \ref{S5}.
In Section \ref{S6}, we prove Theorem ~\ref{main} by reducing to
the model from Sections \ref{S1}--\ref{S5} starting from
a $C^r$ diffeomorphism on a manifold,
and applying Hennion's \cite{He} theorem.
For the H\"older spaces, our proof is elementary: it only
uses integration by parts. 
For the Sobolev spaces, we require in
addition a 
standard $L^t$ estimate (Theorem \ref{th:Ta}) for (operator-valued) 
pseudodifferential operators 
with $C^\infty$ symbols $P(\xi)$ depending only on $\xi$.

\smallskip

{\bf Comments.}
To study spectral properties of Perron-Frobenius operators (and  Ruelle transfer 
operators, more generally), 
it is primarily important to find appropriate spaces for them to act on.
For $C^r$ expanding dynamical systems ($r>1$) and $C^{r-1}$ weights, Ruelle\cite{Ru1}, and later Fried\cite{Fr1} and
Gundlach-Latushkin\cite{GuLa}, showed that the Banach space of 
$C^{r-1}$ functions worked nicely. 
For Anosov diffeomorphisms usual function spaces 
do not work. 
A remedy for this since the seventies is reduction to the expanding case, by taking 
(at least morally) a quotient along the stable 
foliation. This, however, limits the results severely 
 since the stable foliations are in general only H\"older
even if $r=\infty$. 
In the early nineties, Rugh\cite{Rugh}, and then Fried\cite{Fr2}, introduced some ideas which allow to 
bypass this reduction in the case of {\it analytic} Anosov diffeomorphisms. 
These results, together with the work of Kitaev \cite{Ki}
on the radius of convergence of dynamical Fredholm determinants
for hyperbolic systems with finite differentiability, suggested 
that  appropriate  spaces of distributions could be constructed for
$C^r$ hyperbolic dynamics. 
The first major achievement in this direction was made in the work \cite{BKL} by
Blank, Keller and Liverani, in which they considered a Banach space of distributions and 
gave a bound on the essential spectral radius of the Perron-Frobenius operators acting on it. 
However   their methods 
only allowed to exploit limited smoothness of the diffeomorphisms.
In 2004,  Gou\"ezel and Liverani \cite[v1]{GL} 
improved the argument in \cite{BKL} and introduced a new Banach space  of distributions. 
More recently, in \cite[v2]{GL}, they removed technical assumptions in the first version: 
their results now are similar to ours. 
Also in 2004, the first-named author \cite{Ba} 
gave a prototype of the use of Fourier coordinates that we exploit in the
present paper, 
but under a strong assumption on the dynamical foliations. The present paper
is also partly motivated by the argument in \cite{AGT}.
Finally, note that
the spaces of distributions of Gou\"ezel and Liverani
\cite{GL} are similar in spirit to the dual space of 
our H\"older spaces.
The definition of the function spaces of Gou\"ezel and Liverani looks more geometric than ours. 
Our spaces are natural anisotropic versions of the usual H\"older and Sobolev spaces 
as we will see, and fit better in the standard theory of functional analysis. 

\section{Definition of the anisotropic norms.}
\label{S1}

We recall a few facts on Sobolev  and H\"older norms, which motivate our definition 
of anisotropic norms.

Fix an integer $d \ge 1$ and a $C^\infty$ function $\chi:\real\to [0,1]$ with
\[
\chi(s)=
1, \quad \mbox{for $s\le 1$,}\qquad
\chi(s)=0, \quad \mbox{for $s\ge 2$.}
\] 
For $n\in \mathbb Z_+$, define 
$\chi_n:\real^d\to [0,1]$ as $\chi_n(\xi)=\chi(2^{-n}|\xi|)$
and, setting $\chi_{-1}\equiv 0$,
\begin{align*}
&\psi_n:\real^d\to [0,1],\qquad 
\psi_n(\xi)=\chi_n(\xi)-\chi_{n-1}(\xi) .
\end{align*}
We have  
$1=\sum_{n= 0}^{\infty}\psi_{n}(\xi)$, and   $\supp(\psi_n)\subset 
\{\xi\mid 2^{n-1}\le |\xi|\le  2^{n+1}\}$. 
Also $\psi_n(\xi)=\psi_1(2^{-n+1}\xi)$ for $n\ge 1$. Thus,
for every multi-index $\alpha$, there exists a constant $C_\alpha$ such that  
$\|\partial^\alpha \psi_n\|_{L^\infty} \le C_\alpha 2^{-n|\alpha|}$ for all $n\ge 0$, 
and the inverse 
Fourier transform  of $\psi_n$,
$$\widehat \psi_n(x)=
(2\pi)^{-d} \int  e^{ix\xi} \psi_n(\xi)  d\xi,
$$  
decays rapidly, satisfies $\widehat \psi_n(x)=2^{d(n-1)} \widehat \psi_1(2^{n-1}x)$
for $n\ge1$ and all $x$,
and is bounded uniformly in $n$ with 
respect to the $L^1$-norm.

We decompose each $C^\infty$ function $u:\real^d\to \C$ with compact support as
$u=\sum_{n\ge 0}u_n$ by defining for
integer 
$n\in \mathbb Z_+$, 
\begin{equation}\label{decomp}
u_n(x)=\psi_n(D) u(x):=(2\pi)^{-d}\int e^{i(x-y)\xi}\psi_n(\xi)u(y) dy d\xi
=\widehat \psi_n * u  (x) .
\end{equation}
\begin{remark}
The operator $\psi_n(D)$ in (\ref{decomp}) is 
the ``pseudodifferential operator with symbol $\psi_n$."
We refer 
to the books \cite{Ho} and \cite{Ta2} for more about pseudodifferential 
operators, although our text is self-contained,
except for Theorem~\ref{th:Ta}.   
\end{remark}

{}From now on, we fix a compact subset $K\subset \real^d$ with non-empty interior.
Let  $C^{\infty}(K)$ be the space of complex-valued $C^\infty$ 
functions on $\real^d$ supported 
on $K$.  
For a real number $p$ and
$1< t < \infty$, we define on $C^{\infty}(K)$ the norms
\[
\|u\|\hol{p}=\sup_{n\ge 0} \; 2^{pn}\|u_n\|_{L^\infty}
\quad \mbox{and}\quad \|u\|\sob{p,t}=\left\| \left(\sum_{n\ge 0} 4^{pn}
|u_n |^2\right) ^{1/2}\right\|_{L^t}.
\]
It is known that the norm $\|u\|\hol{p}$ is equivalent to the $C^p$ norm
\[
\|u\|_{C^p}=\max\left\{
\max_{|\alpha|\le [p]}
\sup_{x\in \real^d}
|\partial^\alpha u(x)|,
\max_{|\alpha|= [p]}
\sup_{x\in \real^d}\sup_{y\in \real^d/\{0\}} 
\frac{|\partial^\alpha u(x+y)-\partial^\alpha u(x)|}{|y|^{p-[p]}}
\right\}
\]
provided that $p>0$ is not an integer,
and $\|u\|\sob{p,t}$ is equivalent to the generalized Sobolev norm 
\[
\|u\|_{W^{p,t}}=\left\|(1+\Delta)^{p/2} u\right\|_{L^t}
\]
for any $p\in \real$ and $1<t<\infty$. 
(See \cite[Appendix A]{Ta} for a brief account.) 
The {\it little} H\"older space  $C^{p}_*(K)$ is the completion of $C^{\infty}(K)$ with respect 
to 
the norm $\|\cdot \|\hol{p}$. 
The generalized Sobolev space $W^{p,t}_*(K)$ for $1<t<\infty$ is the completion of $C^{\infty}(K)$ 
with respect to the norm $\|\cdot \|\sob{p,t}$. 
\begin{remark}
The {\it little} H\"older space $C^{p}_*(K)$ for non-integer $p>0$ is the closure of $C^\infty(K)$ with
respect to the $C^p$ norm and is smaller than the Banach
space of $C^p$ functions.  Thus our ``H\"older" terminology is slightly incorrect and the notation $C^{p}_*(K)$ may deviate from the standard usage (cf. \cite{Ta}).
\end{remark}

We are going to introduce  anisotropic versions of the norms and spaces above.   
Let  $\cone_+$ and $\cone_-$ be
closed cones in $\real^d$  with nonempty interiors 
such that $\cone_+\cap\cone_-=\{0\}$. 
Let then
$\varphi_+, \varphi_-:\sphere\to [0,1]$ be
$C^{\infty}$ functions on the unit sphere $\sphere$ in $\real^d$ satisfying
\begin{equation}\label{vp}
\varphi_+(\xi)=
\begin{cases}
1, &\mbox{if $\xi\in \sphere\cap \cone_{+}$;}\\
0, &\mbox{if $\xi\in \sphere\cap \cone_{-}$,}
\end{cases} \qquad 
\varphi_-(\xi)=1-\varphi_+(\xi).
\end{equation}
We shall work with combinations
$\Theta=(\cone_+,\cone_-,\varphi_+,\varphi_-)$ as above.
For another such combination $\Theta'=(\cone'_+, \cone'_-, \varphi'_+, \varphi'_-)$,
we write $\Theta'<\Theta$ if
\[
\mbox{closure} (\real^d\setminus \cone_+)\subset\mbox{interior}
(\cone'_-) \cup \{0\} 
\]
(This implies $\cone'_+\subset \cone_+$ and $\cone'_-\supset \cone_-$ in particular.) 
For $n\in \mathbb{Z_+}$ and $\sigma\in\{+,-\}$,  we define
\[
\psi_{\Theta, n, \sigma}(\xi)=
\begin{cases}
\psi_n(\xi)\varphi_{\sigma}(\xi/|\xi|),\quad&\mbox{ if $n>0$;}\\
\chi_n(\xi)/2,\quad&\mbox{ if $n=0$.}
\end{cases}
\]
Note that the $\psi_{\Theta, n, \sigma}(\xi)$
enjoy similar properties as those of the $\psi_n$,
in particular the $L^1$-norm of the rapidly decaying function
$\widehat \psi_{\Theta, n,\sigma}$ is bounded uniformly in $n$.
For a $C^\infty$ function $u:\real^d\to \C$ with compact support, 
an integer $n\in \mathbb Z_+$, $\sigma\in \{+,-\}$, and 
a combination $\Theta=(\cone_+,\cone_-,\varphi_+,\varphi_-)$, we define
\[
u_{\Theta,n,\sigma}=\psi_{\Theta,n,\sigma}(D) u
=\widehat \psi_{\Theta,n,\sigma} * u.
\]
Since $1=\sum_{n= 0}^{\infty}\sum_{\sigma=\pm}\psi_{\Theta, n, \sigma}(\xi)$ 
by definition, we have $
u=\sum_{n\ge 0}\sum_{\sigma=\pm}u_{\Theta,n,\sigma}$.

Let $p$ and $q$ be real numbers.  For   $u\in C^\infty(K)$, 
we define the anisotropic H\"older norm $\|u \|\hol{\Theta,p,q}$   by
\begin{equation}
\|u\|\hol{\Theta,p,q}=\max\left\{\; 
\sup_{n\ge 0}\; 2^{pn}\|  u_{\Theta,n,+}\|_{L^\infty}, 
\; \;\sup_{n\ge 0} \; 2^{qn}\|u_{\Theta,n,-}\|_{L^\infty}\;\right\},
\end{equation}
and the anisotropic Sobolev norm $\|u \|\sob{\Theta,p,q,t}$ 
for $1<t<\infty$ by
\begin{equation}
\|u\|\sob{\Theta,p,q,t}=\; 
\left\| \left(\sum_{n\ge 0} 4^{pn}|  u_{\Theta,n,+}|^2
+ 4^{qn}|u_{\Theta,n,-}|^2 \right)^{1/2}\right\|_{L^t}.
\end{equation}
Let $C^{\Theta,p,q}_*(K)$ be the completion of $C^{\infty}(K)$ with respect to the norm 
$\|\cdot \|\hol{\Theta, p,q}$. 
Likewise, for  $1<t<\infty$, let $W^{\Theta,p,q,t}_*(K)$  be the completion of 
$C^{\infty}(K)$ with respect to the norm $\|\cdot \|\sob{\Theta, p, q,t}$. 
We will call these spaces $C^{\Theta,p,q}_*(K)$ and $W^{\Theta,p,q,t}_*(K)$ of 
distributions the anisotropic H\"older and Sobolev space respectively. 
\footnote{Note that we relate the spaces to isometric images of 
$L^t$ spaces in Appendix \ref{isometric}.}
In Section~\ref{S6}, we will construct the Banach spaces in Theorem \ref{main} by patching these
H\"older and Sobolev spaces 
using local coordinates.


\section{Preliminaries}\label{prelim}
In studying the anisotropic H\"older and Sobolev norms, it is convenient to work 
in different ``coordinates" that we introduce next.
Let $\Gamma=\{ (n,\sigma)\mid n\in \integer_+ ,\sigma\in\{+,-\}\}$ and put
\[
\C^{\Gamma}=\{(f_{n,\sigma})_{(n,\sigma)\in \Gamma}\mid f_{n,\sigma}\in \C\}.
\]
For real numbers $p$ and $q$, and for $\mathbf{f}=(f_{n,\sigma})_{(n,\sigma)\in \Gamma}$ 
and 
$\mathbf{g}=(g_{n,\sigma})_{(n,\sigma)\in \Gamma}$ in
$\C^\Gamma$, we define a norm associated to
a scalar product
\[
|\mathbf{f}|_{\mathcal{W}^{p,q}}=\sqrt
{(\mathbf{f}, \mathbf{f})_{\mathcal{W}^{p,q}}},\qquad
(\mathbf{f}, \mathbf{g})_{\mathcal{W}^{p,q}}=
\sum_{n=0}^{\infty}
\biggl ( 4^{pn}f_{n,+}\cdot \overline{g_{n,+}}+4^{qn}f_{n,-}\cdot \overline{g_{n,-}} \biggr ),
\]
and a norm
$
|\mathbf{f}|_{\mathcal{C}^{p,q}}=
\max\left\{\sup_{n\ge0} 2^{pn}|f_{n,+}|, \sup_{n\ge 0} 2^{qn}|f_{n,-}|\right\}
$.
We then set
\[
\mathcal{W}^{p,q}=\{\mathbf{f}\in \C^{\Gamma}\mid |\mathbf{f}|_{\mathcal{W}^{p,q}}<\infty\}\quad
\mbox{and}\quad
\mathcal{C}^{p,q}=\{\mathbf{f}\in \C^{\Gamma}\mid |\mathbf{f}|_{\mathcal{C}^{p,q}}<\infty\}.
\]
Recall that $K\subset \real^d$ is a fixed compact set.
The operation
\[
\QQ_{\Theta}u=(\psi_{\Theta, n,\sigma}(D)u)_{(n,\sigma)\in \Gamma}
\]
gives the correspondences 
\[
\QQ_{\Theta}:C^{\Theta,p,q}_*(K)\to L^{\infty}(\real^d, \mathcal{C}^{p,q}),\quad 
\QQ_{\Theta}:W^{\Theta,p,q,t}_*(K)\to L^{t}(\real^d, \mathcal{W}^{p,q}).
\]
If we define norms 
\[
\|\mathbf{u}\|_{p,q,\infty}=\| |\mathbf{u}|_{\mathcal{C}^{p,q}}(x)\|_{L^{\infty}}
\quad\mbox{for $\mathbf{u}\in L^{\infty}(\real^d, \mathcal{C}^{p,q})$}
\]
and
\[
\|\mathbf{u}\|_{p,q,t}=\| |\mathbf{u}|_{\mathcal{W}^{p,q}}(x)\|_{L^{t}}
\quad\mbox{for $\mathbf{u}\in L^{t}(\real^d, \mathcal{W}^{p,q})$},
\]
respectively,
the anisotropic H\"older norm and the Sobolev norms coincide with their 
respective pull-backs  
by $\QQ_{\Theta}$:
\begin{equation}\label{pqt}
\|u\|\hol{\Theta, p,q}=\|\QQ_{\Theta}u\|_{p,q,\infty}
\quad\mbox{and}\quad
\|u\|\sob{\Theta, p,q,t}=\|\QQ_{\Theta}u\|_{p,q,t}.
\end{equation}

\smallskip

The pseudodifferential operator $\psi(D)$ with symbol $\psi\in C^{\infty}_0(\real^d)$ extends
to a continuous operator $\psi(D):L^t(\real^d)\to L^t(\real^d)$ for $1\le t\le \infty$ whose 
operator norm is bounded by $\|\widehat\psi\|_{L^1}$, because 
\begin{equation}\label{Young}
\|\psi(D)u\|_{L^t}=\|\widehat\psi*u\|_{L^t}\le \|\widehat\psi\|_{L^1}\|u\|_{L^t}
\end{equation}
by Young's inequality. We will use the following more general result on operator-valued 
pseudodifferential operators:

\begin{theorem}[{\cite[{Theorem 0.11.F}]{Ta}}]\label{th:Ta}
Let $\mathcal{H}_1$ and $\mathcal{H}_2$ be Hilbert spaces and let  $\mathcal{L}(\mathcal{H}_1,
\mathcal{H}_2)$ be the space of bounded linear operators from $\mathcal{H}_1$ 
to $\mathcal{H}_2$ equipped with the operator norm. 
If $P(\cdot)\in C^{\infty}(\real^d, \mathcal{L}(\mathcal{H}_1,
\mathcal{H}_2))$ satisfies
\begin{equation}\label{eqn:szero}
\|D^\alpha_\xi P(\xi)\|_{\mathcal{L}(\mathcal{H}_1,
\mathcal{H}_2)}\le
C_{\alpha}(1+|\xi|^2)^{-|\alpha|/2}
\end{equation}
for each multi-index $\alpha$, then for each  $1<t<\infty$ the operator
\[
P(D):L^t(\real^d,\mathcal{H}_1)\to L^t(\real^d,\mathcal{H}_2) 
\]
is bounded.
\end{theorem}

\begin{remark}
The operator-valued pseudodifferential operator $P(D)$  is defined by 
\[
P(D)u (x) =(2\pi)^{-d}\int e^{i(x-y)\xi}P(\xi) u(y) d\xi dy. 
\]
\end{remark}
\begin{remark}
The proof of Theorem \ref{th:Ta} does not need much knowledge on the theory of
pseudodifferential operators and, in fact,  is rather simple. 
Since the case $t=2$ is proved by using Parseval's identity, one only has to 
check that the arguments in Sections 0.2 and 0.11 of  \cite{Ta} 
extend straightforwardly to the operator-valued case.  
\end{remark}

By ``integration by parts on $w$," we will mean application, 
for $f\in C^2(\real^d)$ and  $g\in C^1_0(\real^d)$,
with $\sum_{j=1}^{d}(\partial_j f)^2\ne 0$ on $\supp(g)$,
of the formula
\begin{align*}
 \int e^{if(w)}g(w) dw
&=
-\sum_{k=1}^d 
\int i(\partial_k f(w))e^{if(w)}\cdot 
\frac{i(\partial_k f(w))\cdot g(w)}
{\sum_{j=1}^{d}(\partial_j f(w))^2} dw
\\
&=
i\cdot \int e^{if(w)}\cdot \sum_{k=1}^d \partial_k\left(\frac
{\partial_k f(w)\cdot  g(w)}
{\sum_{j=1}^{d}(\partial_j f(w))^2}\right) dw ,
\end{align*}
where $w=(w_k)_{k=1}^{d}\in \real^d$, and $\partial_k$ denotes partial differentiation 
with respect to $w_k$. 

If  $f\in C^{1+\delta}(\real^d)$ and  $g\in C^\delta_0(\real^d)$,
for $\delta\in (0,1)$, and $\sum_{j=1}^{d}(\partial_j f)^2\ne 0$ on $\supp(g)$,
we shall consider the following
``regularised integration by parts:" 
\footnote{We thank S. Gou\"ezel for suggesting this.}
Set, for
$k=1,\ldots, d$
\begin{align*}
h_k:=
\frac{i(\partial_k f(w))\cdot g(w)}
{\sum_{j=1}^{d}(\partial_j f(w))^2} .
\end{align*}
Each $h_k$ belongs to $C^\delta_0(\real^d)$.
Let $h_{k,\epsilon}$, for small $\epsilon >0$, be the convolution of
$h_k$ with $\epsilon^{-d} \upsilon(x/\epsilon)$, where 
the $C^\infty$ function
$\upsilon :\real^d \to [0,1]$ is supported in the unit
ball and satisfies $\int \upsilon(x)dx=1$. 
There is $C$, independent of $f$ and $g$,
so that for each small $\epsilon>0$ and all $k$,
$$
\| \partial_k h_{k,\epsilon} \|_{L^\infty}
\le C \|h_k\|_{C^\delta} \epsilon^{\delta-1} ,
\quad
\|  h_k - h_{k,\epsilon} \|_{L^\infty}
\le C \|h_k\|_{C^\delta}  \epsilon^{\delta} .
$$
Finally, for every real number $\Lambda \ge 1$
\begin{align}\label{regparts}
 \int e^{i\Lambda f(w)}g(w) dw
&=-\sum_{k=1}^d \int i\partial_k f(w) e^{i\Lambda f(w)}\cdot  h_k(w) dw\\
\nonumber &= \int \frac{e^{i\Lambda f(w)}}
{\Lambda} \cdot \sum_{k=1}^d \partial_k h_{k,\epsilon} (w) dw \\
\nonumber &\quad-\sum_{k=1}^d 
\int i\partial_k f(w) e^{i\Lambda f(w)}\cdot 
(h_k(w)-h_{k,\epsilon} (w)) dw
.
\end{align}


\section{A pseudolocal property}

Although the pseudodifferential operators $\psi_{\Theta,n,\sigma}(D)$ are not local operators,
i.e., $u_{\Theta,n,\sigma}=\psi_{\Theta,n,\sigma}(D)u$  does not necessarily
vanish  outside of the support of~$u$, we have the following rapid decay
property, which
will be used in Sections \ref{S2} and \ref{Spartition}: 

\begin{lemma}\label{lm:plp}
For all positive real numbers $b$, $c$,  $\epsilon$ and each
$1< t\le \infty$, there exists a constant $C=C(b,c,\epsilon,t)>0$ such that 
\[
|u_{\Theta,n,\sigma}(x)| \le  \frac{C\sum_{\tau=\pm}\sum_{\ell\ge 0}
2^{-c\max\{n,\ell\}}\|u_{\Theta,\ell,\tau}\|_{L^t}}
{d(x,\supp(u))^{b}} ,
\]
for all $n\ge1$,
all $u\in C^\infty(K)$, and all $x\in \real^d$ satisfying $d(x,\supp(u))>\epsilon$.
\end{lemma}

Note that the numerator of the right hand side above is bounded by $C\|u\|\hol{\Theta,p,q}$ in  the case $t=\infty$, and by
 $C\|u\|\sob{\Theta,p,q,t}$ in the case $1<t<\infty$ provided  $c>-q$. 

\begin{proof}
Choose a $C^\infty$ function $\rho:\real^d\to[0,1]$ supported in the disk 
of radius $\epsilon/4$ centered at the 
origin and  so that $\int \rho(x)  dx =1$.  Fix $u\in C^\infty(K)$.
Let $U(\epsilon)$ be the  $\epsilon$-neighborhood of $\supp(u)$. 
Put $\chi_0(x)=\int \mathbf{1}_{U(\epsilon/4)}(y)
\cdot \rho(x-y) dy$, where $\mathbf{1}_{Z}$ denotes the indicator function of 
a subset 
$Z\subset \real^d$. Then $\chi_0$ is supported in $U(\epsilon/2)$, 
with $0\le \chi_0(x)\le 1$ 
for any $x\in \real^d$, and $\chi_0(x)= 1$ for $x\in \supp(u)$. 
Since
$\|\chi_0\|_{C^c_*}$ is bounded by a constant depending only 
on $c$ and  $\epsilon$, we have
\begin{equation}\label{est1}
\|\psi_j(D)\chi_0\|_{L^{\infty}}\le C(c, \epsilon) 2^{-c j}.
\end{equation}
Furthermore, integrating several times by parts on $\xi$  in
\[
\psi_j(D)\chi_0(y)=(2\pi)^{-d} \int e^{i(y-w)\xi}\psi_j(\xi)\chi_0(w) d\xi dw,
\]
we can see that for any $y\in \real^d$ satisfying $d(y,\supp(\chi_0))\ge \epsilon/4$
\begin{equation}\label{est2}
|\psi_j(D)\chi_0(y)|\le C(b,c,\epsilon)\cdot  2^{-c j}
d(y,\supp(\chi_0))^{-b}.
\end{equation}

We assume $d(x,\supp(u))>\epsilon$ henceforth and estimate
\[
\psi_{\Theta,n,\sigma}(D) u(x)=\psi_{\Theta, n, \sigma}(D)(\chi_0 u)(x)
=
\sum_{(\ell,\tau)\in \Gamma} 
\widehat\psi_{\Theta, n, \sigma}* (\chi_0 u_{\Theta, \ell, \tau})(x).
\] 
By the H\"older inequality, we have
\begin{align}\label{hpe}
|\widehat\psi_{\Theta, n, \sigma}*(\chi_0 u_{\Theta, \ell, \tau})(x)|&
\le \|\mathbf{1}_{U(\epsilon/2)}(\cdot)\cdot \widehat\psi_{\Theta, n, \sigma}(x-\cdot)\|_{L^{t'}}
 \|\chi_0 u_{\Theta, \ell, \tau}\|_{L^{t}}\\
 &\le C(b,c,\epsilon, t')\cdot  2^{-c n}
\cdot d(x,\supp(u))^{-b}\cdot \|u_{\Theta, \ell, \tau}\|_{L^{t}}\notag
\end{align}
for any $n$ and $\ell$, where $t'$ is the conjugate exponent of $t$,
i.e. $t^{-1}+(t')^{-1}=1$.  

Suppose that $\ell \ge n+3$. Then we have 
\[
\psi_{\Theta, n, \sigma}(D) ((\psi_j(D)\chi_0)\cdot 
u_{\Theta, \ell, \tau})=0\quad \mbox{for $j<\ell-2$,}
\]
because $\supp(\psi_{\Theta, n, \sigma})$ does not meet $\supp(\psi_j)+\supp(\psi_{\Theta, \ell, \tau})$ which supports the Fourier 
transform of
$(\psi_j(D)\chi_0)\cdot u_{\Theta, \ell, \tau}$. Thus 
\[
\psi_{\Theta, n, \sigma}(D) (\chi_0 u_{\Theta, \ell, \tau})=
\sum_{j\ge \ell-2}\widehat\psi_{\Theta, n, \sigma}*((\psi_j(D)\chi_0)\cdot 
u_{\Theta, \ell, \tau}).
\]
For each $j\ge \ell-2$ with $\ell\ge n+3$, we can see
from  (\ref{est1}-\ref{est2})  that
\begin{align*}
&|\widehat\psi_{\Theta, n, \sigma}*( (\psi_j(D)\chi_0)\cdot 
u_{\Theta, \ell, \tau})(x)|\\
&\qquad\qquad
\le \|\widehat \psi_{\Theta, n, \sigma}\|_{L^{\infty}}\cdot \|\mathbf{1}_{ \real^d \setminus U(\delta)}
\cdot \psi_j(D)\chi_0\|_{L^{t'}}\cdot \|u_{\Theta, \ell, \tau}\|_{L^t}\\
&\qquad\qquad \qquad +
\|\mathbf{1}_{U(\delta)}\cdot\widehat \psi_{\Theta, n, \sigma}(x-\cdot)\|_{L^{t'}}
\cdot \|\psi_j(D)\chi_0\|_{L^{\infty}}\cdot \|u_{\Theta, \ell, \tau}\|_{L^t}
\\
&\qquad \qquad\le C(b,c,\epsilon,t)\cdot   2^{-c j}\cdot d(x,\supp(u))^{-b} 
\cdot \|u_{\Theta, \ell, \tau}\|_{L^{t}},
\end{align*}
where $\delta=\epsilon/2+d(x,\supp(u))/4$. 
(We decomposed the domain of integration in the convolution into 
$U(\delta)$ and its complement.) 
Hence, if $\ell\ge n+3$, we have
\[
|\psi_{\Theta, n, \sigma}(D) (\chi_0 u_{\Theta, \ell, \tau})(x)|
\le  C(b,c,\epsilon,t)\cdot   2^{-c \ell}d(x,\supp(u))^{-b} 
\cdot \|u_{\Theta, \ell, \tau}\|_{L^{t}}.
\]
With this and (\ref{hpe}) we conclude the proof of the lemma. 
\end{proof}

\section{Compact embeddings}
\label{S2}

Recall that $K\subset \real^d$ is a compact subset with non-empty interior.
If $p'\le p$ and $q'\le q$, we have the obvious continuous inclusions
\begin{equation}\label{eqn:inc}
 C_*^{\Theta,p,q}(K)\subset C_*^{\Theta,p',q'}(K),\quad 
W^{\Theta, p,q,t}_*(K)\subset W^{\Theta, p',q',t}_*(K)\quad\mbox{for $1<t< \infty$.}
\end{equation}
Here we prove:

\begin{proposition}\label{prop:cpt}  
If $p'< p$ and $q'< q$, the inclusions (\ref{eqn:inc}) are compact. 
\end{proposition}
\begin{proof}
Take any sequence $u^{(k)}$, $k\ge 1$, in $C^{\infty}(K)$  such that $
\|u^{(k)}\|\hol{\Theta, p,q}<E$ (respectively $\|u^{(k)}\|\sob{\Theta, p,q,t}<E$)
for some positive constant $E>0$. We show that there exists a subsequence $\{k(j)\}$ 
such that $\{u^{(k(j))}\}$ is a Cauchy sequence 
in the norm 
$\|\cdot \|\hol{\Theta, p', q'}$(respectively $\|\cdot \|\sob{\Theta, p', q',t}$). 
For each $(n,\sigma)\in \Gamma$, the Fourier transform $\hat{u}^{(k)}_{\Theta,n,\sigma}$ of  
$u^{(k)}_{\Theta,n,\sigma}$ is a  $C^{\infty}$ function supported on 
$\{\xi\mid 2^{n-1}\le |\xi|\le  2^{n+1}\}$, 
and  its first order derivatives are bounded uniformly 
for $k\ge 1$ and 
$\xi\in \real^d$ since $(1+|x|)u^{(k)}_{\Theta,n,\sigma}(x)$ are uniformly bounded in $L^1$norm from Lemma \ref{lm:plp}. Hence, 
by Ascoli-Arzel\'a's theorem and by the diagonal argument, we can choose a 
subsequence $\{k(j)\}$ such that the sequences 
$\{\hat{u}^{(k(j))}_{\Theta,n,\sigma}\}_{j=0}^{\infty}$  are all Cauchy sequences 
with respect to the $L^\infty$-norm and so is the sequence 
$\{u^{(k(j))}_{\Theta,n,\sigma}\}_{j=0}^{\infty}$.
This is the subsequence with the  required property. 
Indeed, for given $\epsilon>0$, we can choose an integer $N>0$ so that $\sum_{n>N}
(2^{(q'-q)n}+2^{(p'-p)n})E<\epsilon/2$, 
and then we have
\begin{align*}
&\|u^{(k(j))}-u^{(k(j'))}\|\hol{\Theta,p',q'}\\
&\qquad\le
\epsilon/2+
\sum_{n\le  N}\left(
2^{p'n} 
\left\|{u}^{(k(j))}_{\Theta,n,+}-{u}^{(k(j'))}_{\Theta,n,+}\right\|_{L^\infty}+2^{q'n}
 \left\|{u}^{(k(j))}_{\Theta,n,-}-{u}^{(k(j'))}_{\Theta,n,-}\right\|_{L^\infty}\right),
\end{align*}
(respectively the same inequality with the norms $\|\cdot \|\hol{\Theta,p',q'}$ and 
$\|\cdot \|_{L^\infty}$ replaced by $\|\cdot \|\sob{\Theta,p',q',t}$ and $\|\cdot\|_{L^t}$).
The right hand side is $<\epsilon$ for large enough $j$, $j'$.
\end{proof}

\section{A Lasota-Yorke type inequality}
\label{S5}

Let $r >1$.  Let $K, K'\subset \real^d$ be compact subsets with non-empty interiors, and take 
a compact  neighborhood $W$ 
of $K$.  Let $T:W\to K'$ be a $C^{r}$ diffeomorphism onto
its image.
Let $\g:\real^d\to \C$ be a $C^{r-1}$ function such that $\supp(\g)\subset K$. 
In this section we study the transfer operator on $\real^d$:
\[
L: C^{r-1}(K')\to C^{r-1}(K), 
\qquad L u(x)= \g(x)\cdot u\circ T(x).
\]

For two fixed combinations $\Theta=(\cone_+,\cone_-,\varphi_+, \varphi_-)$ and 
$\Theta'=(\cone'_+,\cone'_-,\varphi'_+, \varphi'_-)$ as in Section \ref{S1},
we make   the following {\it cone-hyperbolicity} assumption on $T$:
\begin{equation}\label{conehyp}
DT_x^{tr}(\real^d \setminus \mbox{interior}\,(\cone_{+}))
\subset \mbox{interior}\,(\cone'_{-})\cup\{0\}\quad \mbox{for all $x\in W$,}
\end{equation}
where $DT_x^{tr}$ denotes the transpose of the derivative of $T$ at $x$. 
We put
\[
\|T\|_+=\sup_{x\in \supp(g)}
\sup_{0\neq DT_x^{tr}(\xi)\notin \cone'_{-}} \frac{\|DT_x^{tr}(\xi)\|}
{\|\xi\|}, \quad 
\|T\|_-=\inf_{x\in \supp(g)}
\inf_{0\neq \xi \notin \cone_{+}} \frac{\|DT_x^{tr}(\xi)\|}{\|\xi\|}.
\]

\begin{theorem}\label{th:main2} Fix $\Theta$ and
$\Theta'$ and assume (\ref{conehyp}).
For any $q<0<p$ such that $p-q<r-1$, 
the operator $L$ extends to continuous operators
\[
L: C_*^{\Theta, p,q}(K')\to C_*^{\Theta',p,q}(K),\qquad
L: W_*^{\Theta, p,q,t}(K')\to W_*^{\Theta',p,q,t}(K)
\]
for $1<t<\infty$. 
Furthermore, for any $0\le p'< p$ and $q'< q$ such that $p-q'<r-1$, we have the following 
Lasota-Yorke type inequalities:\\
{\bf H\"older case:} There exist a constant $C$, that does not depend on $T$
or $\g$, 
and a constant $C(T,\g)$, that may depend on $T$ and
$\g$, such that for any $u\in C_*^{\Theta,p,q}(K)$
\begin{equation*}
\|L u\|\hol{\Theta',p,q}\le C\|\g\|_\infty\cdot 
 \max\{ 
 \|T\|_+^{p},\|T\|_-^{q}
\}\|u\|\hol{\Theta,p,q}+C(T,\g)\|u\|\hol{\Theta,p',q'}.
\end{equation*}
\\
{\bf Sobolev case:} For each $1<t<\infty$, there exist a constant $C(t)$, that 
does not depend on 
$T$ or $\g$, 
and a constant $C(T,\g,t)$, that may depend on $T$
and $\g$, such that for any 
$u\in W_*^{\Theta,p,q,t}(K)$
\begin{equation*}
\|L u\|\sob{\Theta',p,q,t}\le C(t)\|\g\|_\infty\cdot 
\frac{
 \max\{ \|T\|_+^{p},
 \|T\|_-^{q}
\}}
{\inf |\det DT|^ {1/t}}
\|u\|\sob{\Theta,p,q,t}+C(T,\g,t)\|u\|\sob{\Theta,p',q',t} .
\end{equation*}
\end{theorem}

For the proof of Theorem \ref{th:main2}, we need more notation. 
By (\ref{conehyp}) there exists a  closed cone
$\widetilde \cone_+$ such that
$\widetilde \cone_+ \subset \mbox{interior} \, (\cone_+)$
and that
\begin{equation}\label{conehyp2}
DT^{tr}_x(\real^d\setminus \mbox{interior}(\widetilde \cone_{+})) \subset  \, \mbox{interior}(\cone'_{-})\cup \{0\} \quad
\mbox{for all $x\in \supp(g)$.}
\end{equation}
Fix also a closed cone $\widetilde \cone_-\subset  \mbox{interior} \, (\cone_-)$ and let 
$\tilde \varphi_+$, $\tilde\varphi_-:\sphere\to [0,1]$ be
$C^{\infty}$ functions  satisfying
\[
\tilde \varphi_+(\xi)=
\begin{cases}
1, &\mbox{if $\xi\notin \sphere\cap \cone_{-}$;}\\
0, &\mbox{if $\xi\in \sphere\cap \widetilde \cone_{-}$.}
\end{cases}, \qquad 
\tilde\varphi_-(\xi)=
\begin{cases}
0, &\mbox{if $\xi\in \sphere\cap \widetilde\cone_{+}$;}\\
1, &\mbox{if $\xi\notin \sphere\cap \cone_{+}$.}
\end{cases}
\]
Recall the function $\chi$ we fixed in the beginning. Put, for $\ell\ge 1$,  
\[
\tilde\psi_\ell(\xi)=
\chi(2^{-\ell-1}|\xi|)-\chi(2^{-\ell+2}|\xi|),
\]
and then define, for $(\ell,\tau)\in\Gamma$, 
\[
\tilde \psi_{\Theta, \ell, \tau}(\xi)=
\begin{cases}
\tilde{\psi}_\ell(\xi)\tilde \varphi_{\tau}(\xi/|\xi|),& \mbox{ if $\ell\ge 1$;}\\
\chi(2^{-1}|\xi|),& \mbox{ if $\ell=0$.}
\end{cases}
\]
Note that $\tilde \psi_{\Theta,\ell,\tau}(\xi)=1$ 
if $\xi \in \supp (\psi_{\Theta,\ell,\tau})$.

Next, fix a closed cone
$\widetilde \cone'_-\subset  \mbox{interior} \, (\cone'_-)$
and take integers $h_{\max}$, $h_{\min}$,
$h_{\min}^{-}$, and 
$h_{\max}^{+}$  such that for all $x\in W$
\begin{align*}
&2^{h_{\min}+4}\|\xi\|< \|DT^{tr}_{x}(\xi)\|< 2^{h_{\max}-4}\|\xi\|
\qquad\mbox{for any $\xi\in \real^d$,}\\
&\phantom{2^{h_{\min}+4}\|\xi\|<} \|DT^{tr}_{x}(\xi)\|< 2^{h_{\max}^{+}-4}\|\xi\|
\qquad\mbox{if $DT^{tr}_{x}(\xi)\notin \widetilde \cone'_{-}$,}\\
&2^{h_{\min}^{-}+4}\|\xi\|< \|DT^{tr}_{x}(\xi)\|
\phantom{< 2^{h_{\max}-4}\|\xi\|}
\qquad\mbox{if $\xi\notin \widetilde \cone_{+}$.}
\end{align*}
By modifying the cones $\widetilde \cone_{+}$ 
and $\widetilde \cone'_-$ if necessary, we may and do assume 
\[
2^{h_{\min}^{-}}> 2^{-5}\cdot \|T\|_{-},\qquad 
2^{h_{\max}^{+}}< 2^{+5}\cdot \|T\|_{+}.
\]

We write $(\ell, \tau) \hookrightarrow (n,\sigma)$ if either 
\begin{itemize}
\item $(\tau,\sigma)=(+,+)$ and $n\le \ell+h_{\max}^+$, or 
\item $(\tau,\sigma)=(-,-)$  and 
$\ell+h_{\min}^{-}\le n$, or
\item $(\tau,\sigma)=(+,-)$ and 
($n\ge h_{\min}^-$ or $\ell \ge - h_{\max}^+$).
\end{itemize}
We write $(\ell, \tau) \not\hookrightarrow (n,\sigma)$ otherwise. 
By the definition of $\not\hookrightarrow$ and by (\ref{conehyp2}),
there exists an integer $N(T)>0$ such that, if 
$(\ell, \tau) \not\hookrightarrow (n,\sigma)$ and 
$\max\{n,\ell\}\ge N(T)$, we have
\begin{equation}\label{lowerbd}
d(\supp(\psi_{\Theta',n,\sigma}), 
DT_{x}^{tr}(\supp(\tilde \psi_{\Theta,\ell,\tau})))
\ge 2^{\max\{n,\ell\}-N(T)} 
\quad \mbox{for $x\in \supp(g)$.}
\end{equation}
Indeed, the case $(\tau,\sigma)=(-,+)$ follows from (\ref{conehyp2}).
Taking
$N(T)\ge \max\{3,h^-_{\min}+3\}$,
the case $(-,-)$ follows from (\ref{conehyp2}) and the definition
of $h_{\min}^-$.  The case $(+,+)$ follows from  the definition
of $h_{\max}^+$ if 
$\xi  \in DT_{x}^{tr}(\supp(\tilde \psi_{\Theta,\ell,\tau}))$
is such that $DT^{tr}_{x}(\xi)\notin \widetilde \cone'_{-}$,
and, taking $N(T)\ge \max\{3,-h^+_{\max}+3\}$, from the fact that 
$\widetilde \cone'_-\subset  \mbox{interior} \, (\cone'_-)$
if $DT^{tr}_{x}(\xi)\in \widetilde \cone'_{-}$.
Finally, the case $(+,-)$ does not occur if we take $N(T)\ge 
\max\{h_{\min}^-,-h_{\max}^+\}$. 

\begin{proof}[Proof of Theorem \ref{th:main2}]
For $v:=L u$, we have
\[
v_{\Theta',n,\sigma}=\sum_{(\ell,\tau)\in \Gamma} \psi_{\Theta',n,\sigma}(D) 
L (u_{\Theta, \ell, \tau}).
\]
We define  $\mathbf{S}$ as 
the formal matrix of operators
\[
S_{n,\sigma}^{\ell,\tau}u
=
\begin{cases}
\psi_{\Theta',n,\sigma}(D) L \, u,&
\mbox{ if $(\ell,\tau)\hookrightarrow (n,\sigma)$},\\
\psi_{\Theta',n,\sigma}(D) L\,\tilde{\psi}_{\Theta,\ell,\tau}(D) u
&\mbox{ if $(\ell,\tau)\not\hookrightarrow (n,\sigma)$},
\end{cases}
\]
for $((\ell,\tau), (n,\sigma))\in \Gamma\times \Gamma$. That is, we set
\[
\mathbf{S}\left((u_{\Theta,\ell,\tau})_{(\ell,\tau)\in \Gamma}\right)
=\left(\sum_{(\ell,\tau)\in \Gamma}S_{n,\sigma}^{\ell,\tau}
u_{\Theta,\ell,\tau}\right)_{(n,\sigma)\in \Gamma}.
\]
Since $\tilde{\psi}_{\Theta,\ell,\tau}(D)u_{\Theta,\ell,\tau}=u_{\Theta,\ell,\tau}$, 
 we have the commutative relation 
$\mathbf{S}\circ \QQ_\Theta=\QQ_{\Theta'}\circ L$. 
For the proof of Theorem \ref{th:main2}, it is enough to show  
\[
\|\mathbf{S}\mathbf{u}\|_{p,q,\infty}<C
2^{\max\{p h_{\max}^+, q h_{\min}^-\}}\|g\|_{L^\infty} \|\mathbf{u}\|_{p,q,\infty}
+C(T,g)\|\mathbf{u}\|_{p',q',\infty}
\]
and that, for $1<t<\infty$,
\[
\|\mathbf{S}\mathbf{u}\|_{p,q,t}
<\frac{C(t)2^{\max\{p h_{\max}^+, q h_{\min}^-\}}\|g\|_{L^\infty}}
{\inf |\det DT|^{1/t}}\|\mathbf{u}\|_{p,q,t}
+C(T,g,t)\|\mathbf{u}\|_{p',q',t} .
\]
To prove the above inequalities, we split the matrix of operator $\mathbf{S}$ into two parts:
\[
\mathbf{S}_0=(\widetilde S_{n,\sigma}^{\ell,\tau}),
\quad
\widetilde S_{n,\sigma}^{\ell,\tau}
=
\begin{cases}
S_{n,\sigma}^{\ell,\tau},&\quad\mbox{if $(\ell,\tau)\hookrightarrow(n,\sigma)$};\\
0&\quad\mbox{if $(\ell,\tau)\not\hookrightarrow (n,\sigma)$},
\end{cases}
\]
and   
$\mathbf{S}_1=\mathbf{S}-\mathbf{S}_0$. 
We first consider $\mathbf{S}_0$. 
This is the composition $\Phi(D)\circ \Psi\circ \mathbf{L}$ of  
\begin{itemize}
\item the operator 
$\mathbf{L}$ defined by $
\mathbf{L}(\mathbf{u})(x)=\g(x)\cdot \mathbf{u}\circ T(x)$,
\item the operator $\Psi$ defined by
\[
\Psi\left((f_{\ell,\tau})_{(\ell,\tau)\in \Gamma}\right)_{(n,\sigma)}=
\sum_{(\ell,\tau)\hookrightarrow(n,\sigma)}f_{\ell,\tau},
\]
where $\sum_{(\ell,\tau)\hookrightarrow(n,\sigma)}$ is the sum over $(\ell,\tau)\in \Gamma$ 
such that $(\ell,\tau)\hookrightarrow(n,\sigma)$, 
\item the pseudodifferential operator $\Phi(D)$ with  symbol 
$\Phi:\real^d\to \mathcal{L}(\C^\Gamma, \C^\Gamma)$,
\[
\Phi(\xi)\left((f_{\ell,\tau})_{(\ell,\tau)\in \Gamma}\right)_{(n,\sigma)}=
\psi_{\Theta',n,\sigma}(\xi) f_{n,\sigma}.
\]
\end{itemize}
Clearly $\|\mathbf{L}\|_{p,q,\infty}\le \|\g\|_\infty$ and $\|\mathbf{L}\|_{p,q,t}\le 
\|\g\|_\infty\sup |\det DT|^{-1/t}$ for $1<t<\infty$. 
Also we can prove 
\begin{equation}\label{b4}
\|\Psi\|_{p,q,t}\le C2^{\max\{p h_{\max}^+, q h_{\min}^-\}}\quad\mbox{ for $1<t\le \infty$}
\end{equation}
as follows. Set  $c(+)=p$, $c(-)=q$ and observe that there is $C$ so that 
\begin{align}\label{cnorm}
&\sum_{(\ell,\tau):
(\ell,\tau)\hookrightarrow (n,\sigma)} 2^{c(\sigma)n-c(\tau)\ell}
\le C2^{\max\{p h_{\max}^+, q h_{\min}^-\}},\, \forall (n,\sigma),\\
\nonumber
& \sum_{(n,\sigma):
(\ell,\tau)\hookrightarrow (n,\sigma)} 2^{c(\sigma)n-c(\tau)\ell}
\le C2^{\max\{p h_{\max}^+, q h_{\min}^-\}},\, \forall (\ell,\tau).
\end{align}
For $\mathbf{f}(x)=(f_{n,\sigma}(x))_{(n,\sigma)\in \Gamma}$, we have, at each point $x\in \real^d$, 
\begin{align}
\label{b2}
&|\Psi(\mathbf{f})|_{C^{p,q}}(x)\le 
C2^{\max\{p h_{\max}^+, q h_{\min}^-\}}
|\mathbf{f}|_{C^{p,q}}(x) 
\intertext
{
and also,  
}
\label{b3}
&|\Psi(\mathbf{f})|_{\mathcal{W}^{p,q}}(x)\le 
C2^{\max\{p h_{\max}^+, q h_{\min}^-\}}
  |\mathbf{f}|_{\mathcal{W}^{p,q}} (x).
\end{align}
The latter inequality is obtained by using
Cauchy-Schwartz and (\ref{cnorm}), as follows:
\begin{align*}
&|\Psi(\mathbf{f})|_{\mathcal{W}^{p,q}}(x)^2
=
\sum_{(n,\sigma)\in \Gamma} 
\left(
\sum_{(\ell,\tau):(\ell,\tau)\hookrightarrow (n,\sigma)}
2^{c(\sigma)n-c(\tau) \ell} 2^{c(\tau) \ell} f_{\ell,\tau}
\right)^2\\
&\le 
\sum_{(n,\sigma)\in \Gamma} 
\bigl(
\sum_{(\ell,\tau):(\ell,\tau)\hookrightarrow (n,\sigma)}
2^{c(\sigma)n-c(\tau) \ell}\bigr)
\bigl( \sum_{(\ell,\tau)\hookrightarrow (n,\sigma)}
2^{c(\sigma)n-c(\tau) \ell} \cdot 2^{2c(\tau) \ell} |f_{\ell,\tau}(x)|^2
\bigr)\\
&\le
C2^{\max\{p h_{\max}^+, q h_{\min}^-\}}
\sum_{(n,\sigma)\in \Gamma} 
 \sum_{(\ell,\tau):(\ell,\tau)\hookrightarrow (n,\sigma)}
2^{c(\sigma)n-c(\tau) \ell} \cdot 2^{2c(\tau) \ell} |f_{\ell,\tau}(x)|^2\\
&\le
C2^{\max\{p h_{\max}^+, q h_{\min}^-\}}
\sum_{(\ell,\tau)\in \Gamma} 
 \sum_{(n,\sigma):(\ell,\tau)\hookrightarrow(n,\sigma)}
2^{c(\sigma)n-c(\tau) \ell} \cdot 2^{2c(\tau) \ell} |f_{\ell,\tau}(x)|^2 \\
&\le C 2^{\max\{p h_{\max}^+, q h_{\min}^-\}}
\sum_{(\ell,\tau)\in \Gamma}  2^{2c(\tau) \ell} |f_{\ell,\tau}(x)|^2 .
\end{align*}

Taking the supremum and $L^t$ norm of  both sides of
(\ref{b2}), (\ref{b3}), respectively, we obtain (\ref{b4}). 
The operator $\Phi(D)$ is bounded with respect to the norm $\|\cdot\|_{p,q,t}$ for 
$1<t\le \infty$: 
If $t=\infty$, this follows from (\ref{Young}) since 
$\widehat\psi_{\Theta',n,\sigma}$ is bounded uniformly for $(n,\sigma)\in \Gamma$ in $L^1$-norm, and
the case $1<t<\infty$ follows from Theorem \ref{th:Ta}. 
Thus we conclude 
\[
\|\mathbf{S}_0\|_{p,q,\infty}\le 
C\|\g\|_\infty\cdot 2^{\max\{p h_{\max}^+, q h_{\min}^-\}},
\]
and
\[
\|\mathbf{S}_0\|_{p,q,t}\le 
\frac{C(t)\|\g\|_\infty\cdot 2^{\max\{p h_{\max}^+, q h_{\min}^-\}}}
{\inf |\det DT|^{1/t}}\quad \mbox{for $1<t<\infty$. }
\]

Next we consider  $\mathbf{S}_1$. 
It only remains to show the following two estimates:
\[
\|\mathbf{S}_1\mathbf{u}\|_{p,q,\infty}<C(T,g)\|\mathbf{u}\|_{p',q',\infty}
\,\, \mbox{
and}\,\, \,
\|\mathbf{S}_1\mathbf{u}\|_{p,q,t}<C(T,g,t)\|\mathbf{u}\|_{p',q',t}, 
\forall 1 < t < \infty .
\]
For this, it is enough to prove that for $1<t\le \infty$,
\begin{equation}\label{eqn:s}
\|S_{n,\sigma}^{\ell,\tau}u\|_{L^t}\le C(T,g)2^{-(r-1)\max\{n,\ell\}}\|u\|_{L^t}\quad 
\mbox{if $(\ell,\tau)\not\hookrightarrow(n,\sigma)$.}
\end{equation}
Indeed, setting $c'(+)=p'$, and $c'(-)=q'$, (\ref{eqn:s}) implies that  
\begin{align*}
\|\mathbf{S}_1\mathbf{u}\|_{p,q,\infty}&
\le \sup_{(n,\sigma)\in \Gamma}
\sum_{(\ell,\tau)\not\hookrightarrow(n,\sigma)}
2^{c(\sigma) n}\|S_{n,\sigma}^{\ell,\tau} u_{\Theta,\ell,\tau}\|_{L^\infty}
\\
&\le C(T,g) \cdot 
\sup_{(n,\sigma)\in \Gamma} \left(\sum_{(\ell,\tau)\not\hookrightarrow(n,\sigma)}
 2^{c(\sigma) n-c'(\tau) \ell-(r-1)\max\{n,\ell\}}\right)\|\mathbf{u}\|_{p',q',\infty},
\end{align*}
and, for $1<t<\infty$,   (in the first inequality below, the triangle inequality is used twice, pointwise and for $L^t$)
\begin{align*}
\|\mathbf{S}_1\mathbf{u}\|_{p,q,t}&
\le \sum_{(n,\sigma)\in \Gamma}
\sum_{(\ell,\tau)\not\hookrightarrow(n,\sigma)}
2^{c(\sigma) n}\|S_{n,\sigma}^{\ell,\tau} u_{\Theta,\ell,\tau}\|_{L^t}
\\
&\le C(T,g) \cdot 
\left(\sum_{(n,\sigma)\in \Gamma} \sum_{(\ell,\tau)\not\hookrightarrow(n,\sigma)}
 2^{c(\sigma)n-c'(\tau) \ell-(r-1)\max\{n,\ell\}}\right)\|\mathbf{u}\|_{p',q',t}.
\end{align*}
The sums in $(\cdot)$ above are finite from the assumption $p-q'<r-1$.

We prove (\ref{eqn:s}). Since (\ref{eqn:s}) is obvious when $\max\{n,\ell\}< N(T)$, we will assume $\max\{n,\ell\}\ge N(T)$.
Rewrite the operator $S_{n,\sigma}^{\ell,\tau}$ in the case 
$(\ell,\tau)\not\hookrightarrow(n,\sigma)$ as 
\[
(S_{n,\sigma}^{\ell,\tau}u)(x)
=(2\pi)^{-2d}\int V_{n,\sigma}^{\ell,\tau}(x,y) \cdot u\circ T(y) |\det DT(y)|
dy,
\]
where 
\begin{equation}\label{Vkernel}
V_{n,\sigma}^{\ell,\tau}(x,y)=\int e^{i(x-w)\xi+i(T(w)-T(y))\eta} 
\g(w)\psi_{\Theta',n,\sigma}(\xi)
\tilde{\psi}_{\Theta,\ell,\tau}(\eta)dw d\xi d\eta.
\end{equation}
Since $\| u\circ T \cdot |\det DT|\|_{L^t}\le C(T,g)\|u\|_{L^t}$, the inequality (\ref{eqn:s}) 
follows if we show that  there exists
$C(T)$ such that for all $(\ell,\tau)\not\hookrightarrow (n,\sigma)$
and all $1<t\le \infty$
the operator norm of the integral operator 
\[
H^{\ell,\tau}_{n,\sigma}:v\mapsto \int V_{n,\sigma}^{\ell,\tau}(\cdot,y)v(y) dy
\]
acting on $L^t(\real^d)$ is 
bounded by $C(T,g)\cdot 2^{-(r-1)\max\{n,\ell\}}$. 

Define the positive-valued integrable function $b:\real^d\to \real$ by
\begin{equation}\label{convol}
b(x)=
1\quad\mbox{ if $|x|\le 1$}, \qquad
b(x)=|x|^{-d-1}\quad\mbox{ if $|x|> 1$.}
\end{equation}
The required estimate on $H^{\ell,\tau}_{n,\sigma}$ follows if we show
\begin{equation}\label{eqn:Kernelest}
|V_{n,\sigma}^{\ell,\tau}(x,y)|\le C(T,g) 2^{-(r-1)\max\{n,\ell\}}\cdot  
2^{d\min\{n,\ell\}}b(2^{\min\{n,\ell\}}(x-y))
\end{equation}
for some  $C(T,g)>0$ and all  $(\ell,\tau)\not\hookrightarrow (n,\sigma)$.
Indeed, as the right hand side 
of (\ref{eqn:Kernelest}) is written as a function of $x-y$, say $B(x-y)$, we have, 
by Young's inequality,   
\[
\|H^{\ell,\tau}_{n,\sigma}v\|_{L^t}\le \|B * v\|_{L^t}\le  \|B\|_{L^1} \|v\|_{L^t}
\le C(T) 2^{-(r-1)\max\{n,\ell\}}\cdot \|b\|_{L^1}\cdot \|v\|_{L^t}.
\]

Below we prove the estimate (\ref{eqn:Kernelest}). 
Integrating (\ref{Vkernel}) by parts $[r]-1$ times on $w$
(in particular, if $1<r<2$ we do nothing), we obtain 
\begin{equation}\label{firststep}
V_{n,\sigma}^{\ell,\tau}(x,y)=
\int e^{i(x-w)\xi+i(T(w)-T(y))\eta} 
F(\xi,\eta,w)\psi_{\Theta',n,\sigma}(\xi)
\tilde{\psi}_{\Theta,\ell,\tau}(\eta)dwd\xi d\eta,
\end{equation}
where $F(\xi,\eta,w)$ is a $C^{r-[r]}$ function in $w$
which is $C^\infty$ in the variables  $\xi$ and $\eta$.
Using (\ref{lowerbd}), we can see that for all $\alpha$, $\beta$
\begin{equation}\label{exest}
\|\partial_\xi^\alpha\partial_\eta^\beta
F(\xi,\eta,\cdot)\|_{C^{r-[r]}}\le C_{\alpha,\beta}(T,g)
2^{-n|\alpha|-\ell|\beta|-([r]-1)\max\{n,\ell\}}.
\end{equation}

Assume first that $r$ is an integer (then, $r=[r]\ge 2$).
Put 
$$G_{n,\ell}(\xi,\eta,w)=F(\xi,\eta,w)\psi_{\Theta',n,\sigma}(\xi)
\tilde{\psi}_{\Theta,\ell,\tau}(\eta).
$$
Consider the scaling
$$
\widetilde G_{n,\ell}(\xi,\eta,w)=G_{n,\ell}(2^{n-1}\xi, 2^{\ell-1}\eta,w).
$$
Then, denoting by $\FF$ the inverse Fourier transform with respect to
the variable $(\xi,\eta)$, we have
$$
\FF G_{n,\ell}(u,v,w)=2^{(n-1)d+(\ell-1)d}\FF 
\tilde G_{n,\ell}(2^{n-1}u,
2^{\ell-1}v,w).
$$
The estimate (\ref{exest}) implies that
$$
\|\partial_\xi^\alpha\partial_\eta^\beta
\widetilde G_{n,\ell}(\xi,\eta,\cdot)\|_{C^{r-[r]}}\le C_{\alpha,\beta}(T,g)
 2^{-([r]-1)\max\{n,\ell\}}\, .
$$
In other words, 
the functions $2^{([r]-1)\max\{n,\ell\}}\FF \widetilde
G_{n,\ell}(u,v,w)$
as functions of $u$ and $v$ are smooth rapidly decaying
functions, uniformly bounded (as rapidly decaying functions) with respect
to $w$, $n$, and $\ell$. By definition, 
\begin{align}
&|V_{n,\sigma}^{\ell,\tau}(x,y)|\le \int |\FF G_{n,\ell}(x-w,T(w)-T(y),w) |dw\\
\nonumber &\quad \le C \int 2^{(n-1)d +(\ell-1)d}
|\FF \widetilde G_{n\ell}(2^{n-1}(x-w), 2^{\ell-1}(T(w)-T(y)),w)| dw  .
\end{align}
With this, it is not difficult to conclude  (\ref{eqn:Kernelest}) 
for integer $r\ge 2$. 

If $r>1$ is not an integer, we 
start from (\ref{firststep}) and rewrite $V_{n,\sigma}^{\ell,\tau}(x,y)$
as  
\begin{equation}\label{firststep'}
\int e^{i\Lambda(x-w)(\xi/\Lambda)+i\Lambda(T(w)-T(y))(\eta/\Lambda)} 
F(\xi,\eta,w)\psi_{\Theta',n,\sigma}(\xi)
\tilde{\psi}_{\Theta,\ell,\tau}(\eta)dwd\xi d\eta,
\end{equation}
for
$\Lambda=2^{\max\{\ell,n\}}$. Recalling (\ref{regparts}),
we apply to (\ref{firststep'})  one regularised integration by parts
for $\delta=r-[r]$ (noting that $T$ is $C^{1+\delta}$).
We get two terms $F_{1,\epsilon}(\xi,\eta,w)$
and $F_{2,\epsilon}(\xi,\eta,w)$.
Choosing   $\epsilon=\Lambda^{-1}$, we may apply the above procedure to each
of them.
\end{proof}

\section{Partitions of unity}
\label{Spartition}

Let $r>0$ and recall $K\subset \real^d$ is compact with nonempty interior.
A $C^r$ partition of unity on $K$
is by definition  a finite family of $C^r$ functions $g_{i}:\real^d\to [0,1]$, $1\le i\le I$, 
such that $\sum_{i} g_i(x)= 1$ for  $x\in K$ and  
$\sum_{i} g_i(x)\le  1$ for  $x\in \real^d$.
The intersection multiplicity of a partition of unity is
$\nu:=\sup_x\#\{i \mid x\in \supp(g_i)\}$.
For $u\in C^\infty(K)$, we set $u_i:=g_i u$ so that $u=\sum_i u_i$. 
In this section, we compare the norms of $u$ and those of the $u_i$'s. 
(This will be useful to refine partitions in the proof of Theorem~\ref{main} in the
next section.)

\begin{lemma} \label{lm:pu}
Let $q\le 0\le p$ satisfy $p-q<r$, and let $p'$ and $q'$ 
be real numbers with $p'<p$
and $q'<q$. 
For every $C^r$ partition of unity $\{g_i\}$ whose intersection multiplicity  is $\nu$, there are constants
$C(\{g_i\})$ and $C(\{g_i\},t)$ (that may depend on the $g_i$'s)
so that for any $u\in C^\infty(K)$
\[
\|u\|\hol{\Theta,p,q}\le \nu\cdot  
\max_{1\le i\le I} \|u_i\|\hol{\Theta,p,q}+C(\{g_i\}) \sum_{1\le i\le I} \|u_i\|\hol{\Theta,p',q'}
\]
and, for all  $1<t<\infty$
\[
\|u\|\sob{\Theta,p,q,t}\le \nu \cdot
\left[ \sum_{1\le i\le I} \|u_i\|\sob{\Theta,p,q,t}^t\right]^{1/t}
+C(\{g_i\},t)
\sum_{1\le i\le I} \|u_i\|\sob{\Theta,p',q',t} .
\]
\end{lemma}

\begin{proof} Let $U(i,\epsilon)$ be the $\epsilon$-neighborhood of the support of $g_i$.
Take $\epsilon>0$ so small that 
the intersection multiplicity of the sets $U(i,\epsilon)$ is $\nu$. 
Decompose $\QQ_{\Theta} u_i$ (recall Section \ref{prelim}) into 
\[
\mathbf{u}^{\mathrm{body}}_i=\mathbf{1}_{U(i,\epsilon)}\cdot \QQ_{\Theta}u_i\quad \mbox{and}\quad
\mathbf{u}^{\mathrm{tail}}_i=\QQ_{\Theta}u_i-\mathbf{u}^{\mathrm{body}}_i.
\]
On the one hand, Lemma \ref{lm:plp} implies
\[
\|\mathbf{u}^{\mathrm{tail}}_i\|_{p,q,\infty}\le C \|u_i\|\hol{\Theta,p',q'}\quad \mbox{and}\quad
\|\mathbf{u}^{\mathrm{tail}}_i\|_{p,q,t}\le C(t) \|u_i\|\sob{\Theta,p',q',t}.
\]
On the other hand, since the intersection multiplicity is $\nu$, we have
\[
\left\|\sum_{i}\mathbf{u}_i^{\mathrm{body}}\right\|_{p,q,\infty}\le 
\nu \cdot \max_i \left\|\mathbf{u}_i^{\mathrm{body}}\right\|_{p,q,\infty},
\]
and, using the H\"older inequality,
\[
\left\|\sum_{i}\mathbf{u}_i^{\mathrm{body}}\right\|_{p,q,t}\le 
\nu^{1/t'} \cdot 
\left[ \sum_i \left\|\mathbf{u}_i^{\mathrm{body}}\right\|_{p,q,t}^t\right]^{1/t}.
\]
Therefore we obtain the estimates in the lemma by using (\ref{pqt}).
\end{proof}

The next proposition gives bounds in the opposite direction. 
\begin{proposition}\label{prop:pu}
Let $q\le 0\le p$, and let $p'$ and $q'$ 
be real numbers with $p'<p$, $q'<q$ and $p-q'<r$. 
If $\Theta'< \Theta$, there are constants $C_0$ and $C_0(t)$
so that   for every $C^r$ partition of unity $\{g_i\}$  there are
constants $C(\{g_i\})$ and $C(\{g_i\},t)$ (which may depend 
on the $g_i$'s) so that for all $u\in C^{\infty}(K)$
\[
\max_{1\le i\le I} \|u_i\|\hol{\Theta',p,q}\le
 C_0\|u\|\hol{\Theta,p,q} +C(\{g_i\})\|u\|\hol{\Theta,p',q'},
\]
and, for $1<t<\infty$, 
\[
\left[\sum_{1\le i\le I}\|u_i\|\sob{\Theta',p,q,t}^{t}\right]^{1/t}
\le 
C_0(t)\|u\|\sob{\Theta,p,q,t} +
C(\{g_i\},t)\|u\|\sob{\Theta,p',q',t} .
\]
\end{proposition}

\begin{proof}
We revisit the proof of Theorem~\ref{th:main2}, setting $T=id$. (Note that assumption (\ref{conehyp}) 
holds since we are assuming $\Theta'<\Theta$.)
Recall $\Phi(D)$ and $\Psi$  there, and let $\mathbf{S}^{(i)}$, $\mathbf{S}^{(i)}_0$, $\mathbf{S}^{(i)}_1$ and $\mathbf{L}^{(i)}$ be  the operators defined  in the same way as $\mathbf{S}$, $\mathbf{S}_0$, $\mathbf{S}_1$ and $\mathbf L$ respectively with $g$ replaced by $g_i$.  
Obviously $|\mathbf{L}_i(\mathbf {f})|_{C^{p,q}}(x)\le |\g_i(x)|
|\mathbf{f}|_{C^{p,q}}(x)$ and $ |\mathbf{L}_i(\mathbf{f})|_{\mathcal{W}^{p,q}}(x)\le 
|\g_i(x)|  |\mathbf{f}|_{\mathcal{W}^{p,q}} (x)$ at each point $x$.
These and (\ref{b2}-\ref{b3}) imply  
\[
\max_{i} \|\Psi\circ \mathbf{L}^{(i)}(\mathbf{u})
\|_{p,q,\infty}\le C_1 \|\mathbf{u}\|_{p,q}\]
and 
\[
\left[\sum_i \|\Psi\circ \mathbf{L}^{(i)}
(\mathbf{u})\|_{p,q,t}^t\right]^{1/t}\le C_1
\|\mathbf{u}\|_{p,q,t}, \forall 1<t<\infty
\]
for $\mathbf{u}=(u_{\Theta,n,\sigma})_{(n,\sigma)\in \Gamma}$. 
By boundedness of $\Phi(D)$, the same estimates hold with $\Psi\circ \mathbf{L}^{(i)}$ replaced by  $\mathbf{S}^{(i)}_0=\Phi(D)\circ \Psi\circ \mathbf{L}^{(i)}$. The conclusion of the proposition then follows from those estimates and the estimates on the operators $\mathbf{S}_1^{(i)}$ parallel to that on $\mathbf{S}_0$ 
in the proof of Theorem~\ref{th:main2}.
\end{proof}


\section{Transfer operators for hyperbolic diffeomorphisms}
\label{S6}

In this section we prove
Theorem~\ref{main} by reducing  to the model of
Sections \ref{S1}--\ref{Spartition}.

\begin{proof}[Proof of  Theorem~\ref{main}]
We first define the  spaces  $C_*^{p,q}(T,V)$ and $W_*^{p,q,t}(T,V)$, 
by using local charts
to patch the  anisotropic H\"older and Sobolev spaces from
Section ~\ref{S1}. Fix 
a finite system of $C^{\infty}$ local charts $\{(V_j, \kappa_j)\}_{j=1}^{J}$ 
that cover the compact isolating neighborhood  $V$ of $\Omega$, 
and a finite system of
pairs of closed cones 
\footnote{We regard $\cone_{j,\pm}$ as constant cone fields in the {\em cotangent} bundle 
$T^*\real^d$.}  
$\{(\cone_{j,+}, \cone_{j,-})\}_{j=1}^{J}$ in $\real^d$
with the 
properties that for all  $1\le j,k\le J$:

\renewcommand{\labelenumi}{(\alph{enumi})}
\begin{enumerate}
\item The closure of $\kappa_j(V_j)$ is a compact subset of $\real^d$.
\item The cones $\cone_{j,\pm}$ are transversal to each other: 
$\cone_{j,+} \cap \cone_{j,-}=\{0\}$.
\item If $x\in V_j\cap \Omega$, the cones 
$(D\kappa_j)^{*}(\cone_{j,+})$ and 
$(D\kappa_j)^{*}(\cone_{j,-})$ in the cotangent space 
contain the normal subspaces of $E^s(x)$ and $E^u(x)$, respectively. 
\item  If 
$T^{-1}(V_k)\cap V_j \neq \emptyset$, setting 
$U_{jk}=\kappa_j(T^{-1}(V_k)\cap V_j)$,
the map in charts 
\[
T_{jk}:=\kappa_k\circ T\circ \kappa_j^{-1}:U_{jk} \to \real^d
\]
enjoys the cone-hyperbolicity condition:
\begin{equation}\label{charthyp}
DT_{jk,x}^{tr}(\real^d\setminus \mbox{interior } (\cone_{k,+}))
\subset \mbox{interior} (\cone_{j,-} )\cup\{0\}, \quad
\forall x \in U_{jk} .
\end{equation}
\end{enumerate}

Choose $C^\infty$ functions $\varphi_j^+, \varphi_j^-: \sphere\to [0,1]$ 
for $1\le j\le J$ 
which satisfy condition (\ref{vp}) with $\cone_\pm=\cone_{j,\pm}$, giving  combinations 
$\Theta_j=(\cone_{j,+}, \cone_{j,-}, \varphi_j^+,\varphi_j^-)$ 
as in Section~\ref{S1}. 
Choose finally a $C^{\infty}$ partition of the unity
$\{\phi_j\}$ on $V$ subordinate to the covering $\{V_j\}_{j=1}^J$, that is, the support of 
each $\phi_j:X\to [0,1]$ is contained in $V_j$ and 
we have $\sum_{j=1}^{J} \phi_j\equiv 1$ on $V$. 

We define the Banach spaces $C_*^{p,q}(T,V)$ and $W^{p,q,t}_*(T,V)$ 
for $1<t<\infty$, respectively, to be 
the completion of $C^{\infty}(V)$ for the norm
\[
\|u\|_{C_*^{p,q}(T,V)}:=\max_{1\le j\le J} 
\|(\phi_j\cdot u)\circ \kappa_j^{-1}\|\hol{\Theta_j, p,q}
\]  
and
\begin{equation}\label{normsob}
\|u\|_{W^{p,q,t}_*(T,V)}:= \max_{1\le j\le J}
\|(\phi_j \cdot u)\circ \kappa_j^{-1}\|\sob{\Theta_j, p,q,t} .
\end{equation}
By this definition, we have that $C_*^{p,q}(T,V)$ and $W^{p,q,t}_*(T, V)$ contain 
$C^s(V)$ for $s>p$ and $W^{p,t}(V)$, respectively, as  dense subsets. 
Take and fix real numbers $0\le p'<p$ and $q'<q$ such that $p-q'<r-1$.  
By Lemma ~\ref{prop:cpt} and a finite
diagonal argument over $\{1, \ldots, J\}$, 
we can see that the inclusions 
$C^{p,q}_*(T, V) \subset C^{p',q'}_*(T, V)$ and 
$W^{p,q,t}_*(T,V) \subset W^{p',q',t}_*(T, V)$ are compact.

For $m\ge 1$ and $j$, $k$ so that
\[
V_{m, jk}:=T^{-m}(V_k)\cap V_j\cap \left(\cap_{i=0} ^{m}T^{-i}(V)\right) \ne \emptyset,
\]
we may consider the map  in charts 
\[
T^m_{jk}=\kappa_k\circ T^m\circ \kappa_j^{-1}:
\kappa_j(V_{m,jk})\to \real^d.
\]
Note that (\ref{charthyp}) implies that
\begin{equation}\label{charthyp2}
(DT^m_{jk,x})^{tr}(\real^d\setminus \mbox{interior } (\cone_{k,+}))
\subset \mbox{interior} (\cone_{j,-} )\cup\{0\}, \quad
\forall x \in \kappa_j(V_{m,jk}) .
\end{equation}

For $1<t\le \infty$, we 
set
\[
 \Lambda_{m,t}= \max_j \max_k \sup_{x\in \kappa_j(V_{m,jk})}  
 \frac{|\g^{(m)}\circ \kappa_j^{-1}(x)|\cdot
 \max\{ (\|T^m_{jk}\|_+(x))^p, (\|T^m_{jk}\|_-(x))^{q} \}}{|\det DT^{m}_x|^{1/t}}
\]
where 
\[
\|T_{jk}^m\|_+(x)=
\sup\left\{ \frac{\|(DT_{jk}^m)_x^{tr}(\xi)\|}{\|\xi\|}\;;\;
 0\neq (DT_{jk}^m)_x^{tr}(\xi)\notin \cone_{j,-}\right\},
\]
and
\[
\|T_{jk}^m\|_-(x)=
\inf\left\{ \frac{\|(DT_{jk}^m)_x^{tr}(\xi)\|}{\|\xi\|}\;;\;
 0\neq \xi \notin \cone_{k,+}\right\}.
\]
Then a standard argument on hyperbolic sets gives a constant $C(t)>1$ that 
does not depend on $m>0$ such that
\begin{equation}\label{eta}
C(t)^{-1} R^{p,q,t}(T,\g,\Omega,m)\le 
\Lambda_{m,t}\le C(t) R^{p,q,t}(T,\g,\Omega,m).
\end{equation}

The definition of $\Lambda_{m,t}$ involves first taking a
maximum and a product, and then taking the supremum over $x$.
We shall apply  Theorem~\ref{th:main2} in a moment: the
upper bound there corresponds to taking a supremum first.
Since different points in $\kappa_j(V_{m,jk})$ may have very different
itineraries, it is necessary to refine our partition of unity, depending
on $m$. This will not cause
problems since 
we can take arbitrarily fine  finite $C^\infty$ partitions of unity on $\real^d$,
with intersection multiplicities  bounded uniformly by a  
constant  depending only on $d$.   Using such a partition
of unity, we decompose the function
$u_{jk}=(\phi_k(\phi_{j} \circ T^{-m}) \cdot u )\circ \kappa_k^{-1} $
into $u_{jk,i}$ for $1\le i \le I_{jk}$.
Take  combinations $ \Theta'_{k}<\Theta_k$  (close to $\Theta_k$)
so that the iterated
cone-hyperbolicity condition (\ref{charthyp2}) holds with $\Theta_k$ replaced by $\Theta'_{k}$. 
For each $m$, by taking a sufficiently fine partition
of unity,
we can apply Theorem~\ref{th:main2} to obtain, for $1\le i\le I_{jk}$, 
\[
\|g^{(m)}\circ \kappa_j^{-1}\cdot u_{jk,i}\circ T^m_{jk}\|\hol{\Theta_j,p,q}
\le 2\Lambda_{m,\infty}\cdot \|u_{jk, i}\|\hol{\Theta'_k,p,q} +C\|u_{jk,i}\|\hol{\Theta'_k,p',q'}.
\]
Then, using Lemma \ref{lm:pu} and Proposition \ref{prop:pu}, we  get
\begin{align*}
\|g^{(m)}\circ \kappa_j^{-1}\cdot u_{jk}\circ T^m_{jk}\|\hol{\Theta_j,p,q}\le 
C_1\cdot \Lambda_{m,\infty}\cdot 
\|u_{jk}\|\hol{\Theta_k,p,q}+C_1(m)\cdot  
\|u_{jk}\|\hol{\Theta_k,p',q'},
\end{align*}
where $C_1$ is a constant that does
not depend on $m$. 
Thus, using Proposition \ref{prop:pu} again, we obtain the following Lasota-Yorke type inequalities:
\begin{align*}
\|\mathcal{L}_{T,g}^m u\|_{C_*^{p,q}(T,V)}\le C_2 \cdot J\cdot \Lambda_{m,\infty}\cdot 
\|u\|_{C_*^{p,q}(T,V)} 
+C_2(m) \| u\|_{C_*^{p',q'}(T,V)} , \, m \ge 1.
\end{align*}
Likewise, we obtain for $1<t<\infty$
\begin{align*}
\|\mathcal{L}_{T,g}^m u\|_{W_*^{p,q,t}(T,V)}\le C_2(t) \cdot J\cdot \Lambda_{m,t}\cdot 
\|u\|_{W_*^{p,q,t}(T,V)} 
+C_2(m,t) \| u\|_{W_*^{p',q',t}(T,V)}.
\end{align*}
Finally Hennion's theorem \cite{He}
gives the claimed upper bounds 
\[
\liminf_{m \to \infty}
(C(t) \Lambda_{m,t})^{1/m}=
R^{p,q,t}(T,\g,\Omega)
\]
for the essential
spectral radius of $\mathcal{L}_{T,g}$. 
\end{proof}

\begin{remark}\label{MTO}
The proof above applies to  (hyperbolic) mixed transfer operators \cite{Ki}.
\end{remark}

\begin{remark}
Though it is not  explicit in our notation, the definition of the spaces 
$C_*^{p,q}(T,V)$ and $W_*^{p,q,t}(T,V)$ depends on the system of charts 
$\{(V_j,\kappa_j)\}_{j=1}^J$, the set of combinations  
$\{(\cone_{j,+}, \cone_{j,-}, \varphi_{j,+},\varphi_{j,-})\}_{j=1}^{J}$, 
and the partition of unity $\{\phi_j\}_{j=1}^J$. Choosing a different system of local charts, 
a different set of combinations, or a different partition of unity, does not a 
priori give rise to equivalent norms, though Theorem \ref{th:main2} gives relations.  
This is a little unpleasant, but does not
cause  problems. 
\end{remark}

\appendix

\section{Relating the anisotropic Banach spaces with
$L^t(K)$}

\label{isometric}

Recalling $\varphi_\pm$ from Section~\ref{S1},
define for real numbers $p$ and $q$ the symbols
\[
\Psi_{\Theta,p,+}(\xi)=(1+|\xi|^2)^{p/2}\varphi_+(\xi/|\xi|)\quad\mbox{and}
\quad \Psi_{\Theta,q,-}(\xi)=(1+|\xi|^2)^{q/2}\varphi_-(\xi/|\xi|),
\]
and, recalling the compact set $K\subset \real^d$ with nonempty interior,
define norms for $u\in C^\infty(K)$ 
and $1<t<\infty$ by
\begin{equation}
\begin{cases}
\|u\|\sob{\Theta, p, q,t}^{\dagger\dagger}&
=\|\Psi_{\Theta,p,+}(D)u+\Psi_{\Theta,q,-}(D)u\|_{L^t},\\
 \|u\|\sob{\Theta, p, q,t}^{\dagger}&
=\|\Psi_{\Theta,p,+}(D)u\|_{L^t}+\|\Psi_{\Theta,q,-}(D)u\|_{L^t}.
\end{cases}
\end{equation}
Let $\|\cdot  \|\sob{\Theta, p, q,t}$ be the
norm defined in Section~\ref{S1}. We shall prove below (using Theorem~\ref{th:Ta}) that 
for each $1<t<\infty$ there is a constant
$C>0$ so that
\begin{equation}\label{eqn:equinorm}
C^{-1}\|u\|\sob{\Theta, p, q,t}\le \|u \|\sob{\Theta, p, q,t}^{\dagger}\le 
C\|u \|\sob{\Theta, p, q,t}\quad\mbox{ for $u\in C^{\infty}(K)$.}
\end{equation}
Let $W^{\Theta,p,q,t}_\dagger(K)$ be the completion of 
$C^{\infty}(K)$ with respect to $\|\cdot \|\sob{\Theta, p, q,t}^{\dagger}$. Then
$W^{\Theta,p,q,t}_\dagger(K)= W^{\Theta,p,q,t}_*(K)$ by (\ref{eqn:equinorm}).

Let $W^{\Theta,p,q,t}_{\dagger\dagger}(K)$ be the completion of 
$C^{\infty}(K)$ with respect to $\|\cdot \|\sob{\Theta, p, q,t}^{\dagger\dagger}$.
Then
\begin{equation}\label{distr}
W^{\Theta,p,q,t}_{\dagger\dagger}(K)
=\{u \in \SS'(\real^d)\mid \supp(u)\subset K\, , \quad 
\|u \|\sob{\Theta, p, q,t}^{\dagger\dagger} < \infty\} , \forall 1<t<\infty .
\end{equation}
Clearly,
\begin{equation}\label{eqn:equinorm4}
\|u \|\sob{\Theta, p, q,t}^{\dagger\dagger}\le \|u \|\sob{\Theta, p, q,t}^{\dagger}
\quad\mbox{ for $u\in C^{\infty}(K)$.}
\end{equation}
Though we do not know whether the norm $\|\cdot \|\sob{\Theta, p, q,t}^{\dagger\dagger}$ 
is equivalent to 
$\|\cdot \|\sob{\Theta, p, q,t}$, we show below (using
Theorem~\ref{th:Ta}) that, for each $1<t<\infty$, 
if $\Theta'>\Theta$ and $p\ge q$,  
\begin{equation}\label{eqn:equinorm2}
\|u \|\sob{\Theta, p, q,t}^{\dagger}\le C\|u \|\sob{\Theta', p, q,t}^{\dagger\dagger} 
\quad\mbox{ for $u\in C^{\infty}(K)$}
\end{equation}
for some constant $C>0$.
From (\ref{eqn:equinorm}), (\ref{eqn:equinorm4}) and (\ref{eqn:equinorm2}), it follows that,  
if $\Theta'>\Theta$ then
for all $1<t<\infty$  and $p\ge q$, 
then 
$W^{\Theta',p,q,t}_*(K)\subset W^{\Theta,p,q,t}_*(K)$, and, also,
$$
W^{\Theta',p,q,t}_{\dagger\dagger}(K)
\subset W^{\Theta,p,q,t}_*(K)\subset W^{\Theta,p,q,t}_ {\dagger\dagger}(K).
$$ 

Let $W^{\Theta,p,q,t}_ {\dagger\dagger}(\real^d)\subset \SS'(\real^d)$ be  the isometric 
image  of the restriction to $L^t(\real^d)$
of the inverse of the bijective operator
$\Psi_{\Theta,p,+}(D)+\Psi_{\Theta,q,-}(D)
:\SS'(\real^d) \to \SS'(\real^d)$.
For $1<t<\infty$, (\ref{distr}) implies that 
$W^{\Theta,p,q,t}_ {\dagger\dagger}(K)$ is just (isometrically)
$\{u \in W^{\Theta,p,q,t}_ {\dagger\dagger}(\real^d)
\mid \supp(u)\subset K\}$. 
The space
$W^{\Theta,p,q,t}_\dagger(K)$  may be described in a
similar (although not as neat)
way using the injective (non surjective)
operator $(\Psi_{\Theta,p,+}(D),\Psi_{\Theta,q,-}(D)):
\SS'(\real^d) \to \SS'(\real^d)\oplus \SS'(\real^d)$.

\smallskip
To finish,  we  prove (\ref{eqn:equinorm}) 
and (\ref{eqn:equinorm2}). 
For the proof of (\ref{eqn:equinorm}), it is enough to show
\begin{equation}\label{eqn:equinorm5}
C^{-1}\|\Psi_{\Theta,p,+}(D)u\|_{L^t}\le 
\left\|\left(\sum_{n\ge 0}(2^{pn}\psi_{\Theta,n,+}(D)u)^2\right)^{1/2}\right\|_{L^t}
\le C\|\Psi_{\Theta,p,+}(D)u\|_{L^t}
\end{equation}
for some constant $C>0$ and the corresponding claim for $\|\Psi_{\Theta,q,-}(D)u\|_{L^t}$. 
Since the norm $\|\cdot \|\sob{0,t}$ is equivalent to the 
norm $\|\cdot \|_{L^t}$ as  noted in Section~\ref{S1}, and 
since $\psi_{n}(\xi)\Psi_{\Theta,p,+}(\xi)=a^p(\xi)\psi_{\Theta,n,+}(\xi)$ for $n\ge 1$,
with 
$a^p(\xi)=(1+|\xi|^2)^{p/2}$, we have
\begin{align*}
&C^{-1}\|\Psi_{\Theta,p,+}(D)u\|_{L^t}\le \\
&\qquad\qquad\left\|\left(\sum_{n\ge 0}(a^p(D)\psi_{\Theta,n,+}(D)u)^2\right)^{1/2}\right\|_{L^t}
=\|\Psi_{\Theta,p,+}(D) u\|_{W_*^{0,t}}\\
&\qquad\qquad\qquad\qquad\qquad\qquad \le C\|\Psi_{\Theta,p,+}(D)u\|_{L^t}.
\end{align*}
It is easy to see  that there exist  
$P, Q\in C^{\infty}(\real^d, \mathcal{L}(\ell^2,\ell^2))$ 
satisfying  (\ref{eqn:szero}) such that the corresponding pseudodifferential 
operators $P(D)$ and $Q(D)$ acting on
$L(\real^d,\ell^2)$ transform $(2^{pn}\psi_{\Theta,n,+}(D)u)_{n\ge 0}$ to 
$(a^p(D)\psi_{\Theta,n,+}(D)u)_{n\ge 0}$ and the reverse. 
Thus Theorem \ref{th:Ta} gives (\ref{eqn:equinorm5}). The corresponding claim for 
$\|\Psi_{\Theta,q,-}(D)u\|_{L^t}$ is shown in a parallel manner. 

In order to prove (\ref{eqn:equinorm2}), it is enough to show  
\[
\max\{\|\Psi_{\Theta,p,+}(D)u\|_{L^t}, \|\Psi_{\Theta,q,-}(D)u\|_{L^t}\}
\le C\cdot \|\Psi_{\Theta',q,-}(D)u+\Psi_{\Theta',p,+}(D) u\|_{L^t}
\]
for some constant $C$.
Since 
$$\Psi_{\Theta,p,+}(\xi)/(\Psi_{\Theta',p,+}(\xi)+\Psi_{\Theta',q,-}(\xi))
\mbox{ and } 
\Psi_{\Theta,q,-}(\xi)/(\Psi_{\Theta',p,+}(\xi)+\Psi_{\Theta',q,-}(\xi))$$ 
both satisfy 
(\ref{eqn:szero}), this follows from Theorem \ref{th:Ta}.

\bibliographystyle{amsplain}

\end{document}